\documentclass[12pt,a4paper]{article}
\usepackage{hyperref,amsmath,amssymb}
\usepackage[round]{natbib}
\usepackage{enumitem}

\usepackage{tipa}
\usepackage{graphicx}
\usepackage{geometry}
\usepackage{parskip}

\author{Rowena Ball\footnote{Contact: Professor Rowena Ball, email Rowena.Ball@anu.edu.au, Mathematical Sciences Institute,  Australian National University, Ngunnawal and Ngambri Lands, Canberra, ACT 2601} 
} 

\date{Mathematical Sciences Institute,  Australian National University, Canberra, Australia}

\begin{document}

\title{Empirical mathematics in Australian Indigenous  Smoke Telegraphy}

\setlength{\parskip}{6pt}
\maketitle
\titlepage

\begin{abstract}\noindent 

Mathematics curriculums at most universities tend to perpetuate a belief that higher mathematics is historically and culturally European. First Nations and minority students may not see their identities and cultures reflected in the discipline, yet university mathematics educators are keen to diversify and broaden the appeal of their courses.  This article presents an investigation on the mathematics of smoke telegraphy, as a contribution to inlaying cross-cultural mathematical heritage in the curriculum. 
Across Indigenous societies of Australia the technology and practice of smoke telegraphy was developed to a sophisticated level over millennia to fill a need for long-distance communications. Through an original bibliographic and archival analysis, we show that smoke signalling and telegraphy used empirical mathematics of symmetries, frequency coding, and understanding of fluid dynamics.  
We juxtapose this applied mathematical knowledge, within context,  against the timeline of Western understanding and development of these strands of mathematics.

\end{abstract}

\newpage

\section{Introduction \label{section1}}

Smoke telegraphy is a wireless-broadcasting,  long-distance communication system 
that was developed by Aboriginal and Torres Strait Islander peoples of the Australian continent and islands and widely used from ancient times to the 1940s. The technology engaged the empirical sciences of preparing and deploying smoke in  precise forms and distinct colours to code, send, decode, and relay messages with high fidelity. The colour (or chemical) science applied in smoke signalling has been investigated in \cite{Ball:2026}, which may be regarded as a companion article to this one. Mathematical perception and know-how, too, were intrinsic to smoke telegraphy, and in this article I elicit some of this knowledge within context. Beyond attributes such as numerical coding, I highlight the use of  chiral symmetry, which was understood and actively created by Indigenous smoke signallers, and  the use of frequency coding of the visible signals.

A widely-held belief, often deeply internalized, 
is  that mathematics is primarily of European\footnote{I shall use  the descriptors `European', `Western' and  `mainstream' interchangeably, or as appropriate. Whatever the term `higher mathematics' includes, it is always culturally European.}  provenance, thence gifted to the curriculums of the whole world as the highest form of knowledge \citep{Hoyrup:1996,Raju:2001}. This belief is not true, of course. For how could one group/race of humans be intrinsically capable of mathematical thought, or have an inherent ability to do mathematics, but other groups/races not? 
Although recognition that all societies developed mathematical  capabilities and practices is belated, appreciation of that knowledge is growing \citep{D'Ambrosio:2000}, albeit  slowly and patchily.

In this atmosphere of uneasy ambivalence, undergraduate students of mathematics may well wonder why their courses are rich with mathematics that was done in Europe and the UK and the US during the last two or three centuries but contain next to nothing about what, or how, mathematics was done by other societies. Currently this curriculum gap is formidably wide, notwithstanding some four decades of research on culturally diverse mathematics that was scoped in~\cite{Xu:2025}. In `math-ing up' the infill, the focus of a single article  is necessarily narrow. 
Here I distill, from primary historical sources,   some of the mathematical knowledge inherent to the practice of smoke telegraphy within and between Indigenous civilizations of the Australian continent and islands, set alongside analogous or cogent mathematical developments in Europe. 

This article is pitched at Indigenous and non-Indigenous students and educators of university mathematics courses. 
Although tremendous advances have been made in decolonising school mathematics \citep{ATSIMA,Edmonds:2026}, 
Western mathematical framings and epistemology, histories and content, heroes and role models are all-pervasive in university mathematics curriculums and textbooks \citep{Bishop:1990,DecolonisingMaths,Pearce:2003}. 
Indigenous and minority students who disengage from mathematics, regarding it as some sort of secret white-man's business irrelevant, or even inimical, 
to their culture and identity,  may be unnecessarily sidelined from high-paying careers for which the degree courses require at least first-year mathematics. In that case the whole of society is poorer, 
as well as the student and their community. 

\newpage

Yet many mathematics educators at this level increasingly are practising innovative methods of engagement and inclusion, and  would like nothing better than to broaden the appeal of their courses and degrees to a greater diversity of students, including those of Indigenous and minority-culture identities \citep{Leung:2020}.\footnote{The existential threat to mathematics departments of sporadic claims by science departments to the effect that ``we can teach our students in-house all the mathematics they need to know'' is nothing if not a powerful motivator for greater inclusivity.}  Usually, though, they lack resources and are time-constrained. (First-year maths course coordinators  typically are responsible for around 300--600 or more students.) On the students' side, many from all backgrounds are seeking a more inclusive mathematical education.\footnote{The author's  one-semester 3rd-4th year accredited course on non-Western and Indigenous mathematics was vastly over-subscribed.}  

One aim of this article is to present the mathematical aspects of smoke telegraphy as material that can underpin topics such as symmetry and the Fourier transform, which are usually encountered in mainstream first year mathematics major and non-major (or service) courses, studied in more depth in second year, and have practical importance in most of the sciences and engineering. 

\paragraph{Overview of article:} 
Smoke telegraphy is described in section \ref{section2}. 
Selected excerpts from primary sources are reproduced in the Appendix. These excerpts are  \textit{kabari banjeeri} (`very good companions', in old Mithaka dialect) of this study, interesting in their own right, and should be read in full in parallel with and as referenced in the main text. The main results are extracted from them and summarized in section~\ref{results} and discussed further in section \ref{discussion}, where I juxtapose the mathematical knowledge evidently known to and practised by the Indigenous master smoke telegraphers against contemporaneous developments in Western mathematical knowledge. Conclusions  are drawn and limitations acceded in section \ref{conclusion}.  

\section{What is smoke telegraphy? \label{section2}} 
The aphorism `gone up in smoke' implies a loss of information, and that something important has vanished or been destroyed, both literally and metaphorically. But to Aboriginal and Torres Strait Islander peoples, meaningful information could be written \textit{into} smoke streams, received visually by distant recipients, read off, and re-coded and relayed on further as necessary. This was smoke telegraphy, an extraordinarily sophisticated and effective long-distance communications technology that was practised throughout the Australian continent and islands. Some of its infrastructure and tactical uses during the various frontier wars of the 19th century in south-east Queensland have been documented by \cite{Kerkhove:2021}. During that period and the first half of the 20th century, white colonial society was fascinated with smoke telegraphy and there were a great many eyewitness reports in newspapers and magazines, and observations described in talks to learned societies, memoirs and books.  

How could information, data, and messages be coded and transmitted in smoke? Sure, smoke is a signal of fire, but is it not an evanescent, dissipative aerosol that dissolves into the sky, even the furious, menacing emissions from a bushfire? 

\newpage

\begin{figure}[ht]
\centerline{\includegraphics[scale=1]{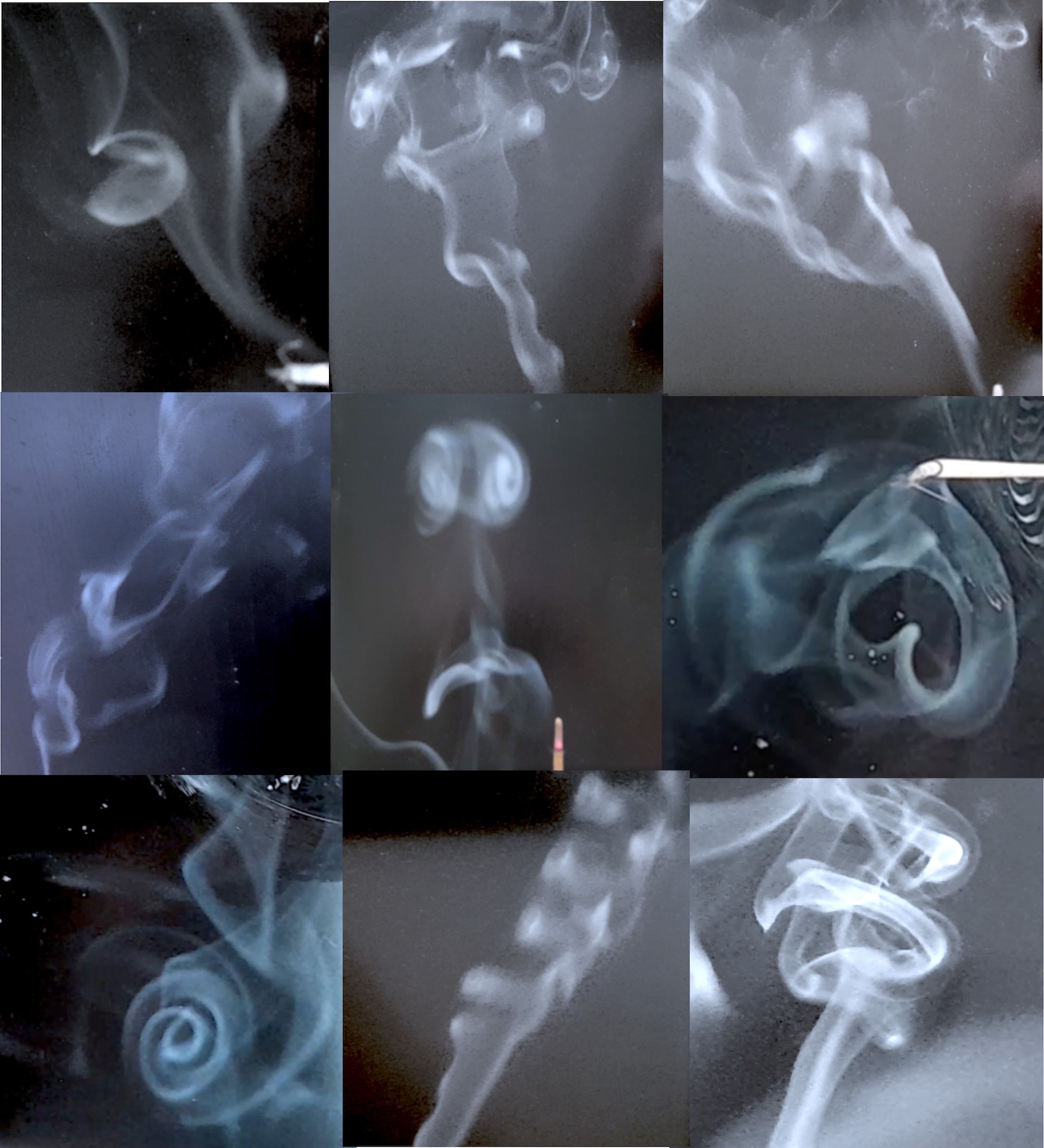}}
\caption{\label{figure1} The vortical nature and natural structure-forming capabilities  of smoke flows are evident in this pastiche of stills from several experimental home movies.}
\end{figure}

Yes and no. A sense of the information coding and carrying potential of smoke as a transmission medium can be gained simply by lighting a stick of incense (figure~\ref{figure1}). In steady conditions you may observe a laminar plume of smoke emanating from the tip, which bifurcates into twin,  counter-rotating vortex streams, due to inevitable perturbations in the interfacial velocity shear at both edges of the smoke flow (aka the Kelvin-Helmholtz instability). These primary vortex streams in turn become unstable, developing into helical, scroll, double-twisted, and chain structures, mushrooms, hearts, pulses and puffs, sheets and Karman streets, thence into all kinds of  pareidolias (see Appendix, excerpt~\ref{a} p. 279) and other intriguing secondary  structures, some of which one may imagine look like an x-ray film of an alien being, until the energy cascade dissipates them at the small scales. (The art-photos of cigarette smoke by professional photographer Thomas Herbrich at  \href{https://www.herbrich.com/photos}{herbrich.com/photos} are also worth viewing.)  

\newpage
\begin{figure}[ht]
\centerline{
\includegraphics[scale=1]{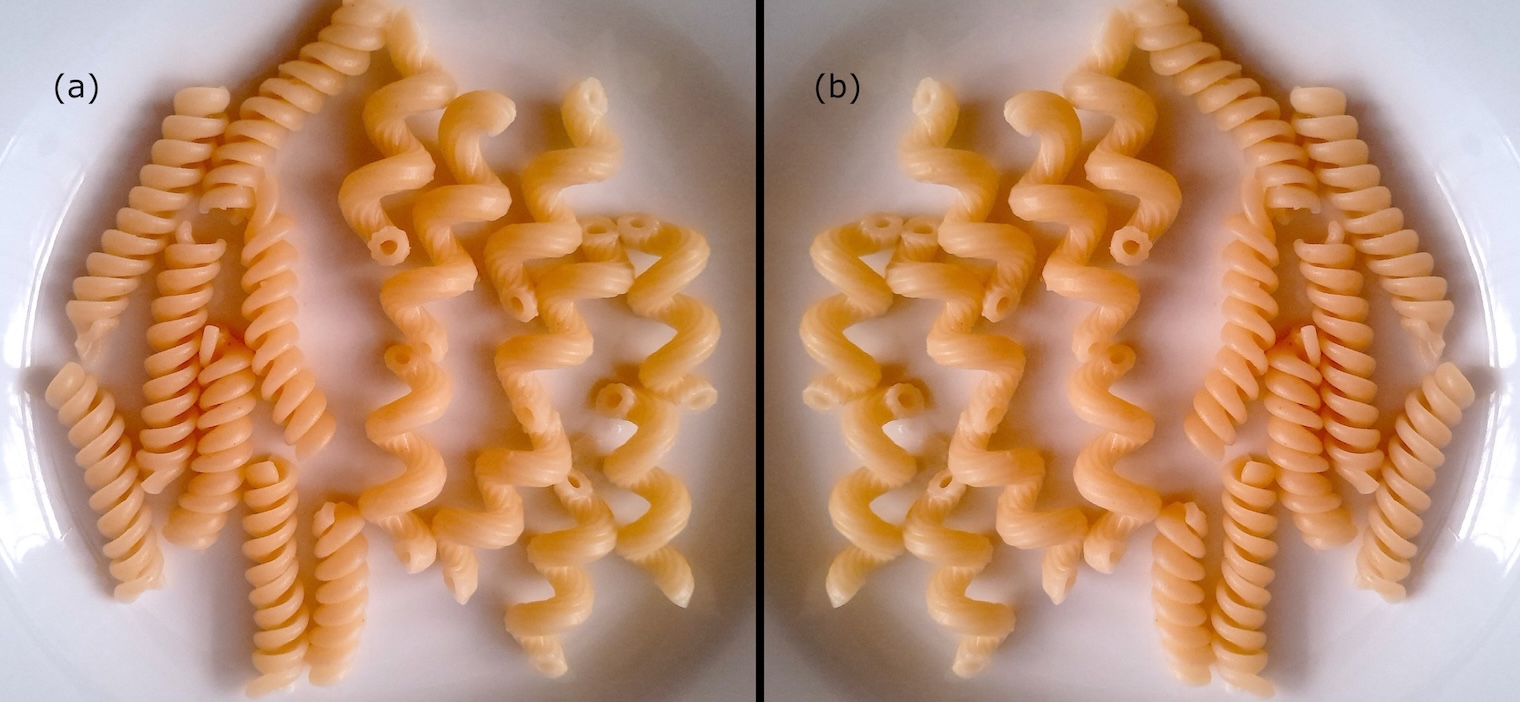}
}
\caption{\label{figure2} Geometry of smoke streams in terms of pasta (ever a cheap and tasty mathematical prop \citep{Hildebrand:2010}). (a) The fusilli shape on the left of the plate is analogous to a primary vortex stream and the  cavatappi  shape on the right of the plate is analogous to the secondary smoke spirals. (b) The chiral asymmetry property is apparent in the mirror image of the plate of pasta---the shapes are non-superimposable on their enantiomorphs in (a).}
\end{figure}
The primary vortices are  shed as left-handed and right-handed pairs. Their 3-dimensional structure embodies the geometric symmetry property of chirality, or handedness  \citep{Darvas:2007}. (Despite being fundamental to biology, this property can be uncomfortable in physical and social terms if misunderstood: try shaking your colleague's left hand with your right hand.) An object is chiral if it cannot be superposed on its real or virtual mirror image. Such an object and its same-but-different mirror self are called a pair of enantiomorphs, or enantiomers if they are molecules  \citep{Adawy:2022}. Smoke vortices and spirals embody axial or helical chirality. 
As the primary vortex streams themselves wobble into secondary helices, their handedness  is passed on too (figure~\ref{figure2}).  

The smoke stream structures in figure \ref{figure1} are generated naturally. Aboriginal and Torres Strait Islander signallers were expert in exploiting the information coding and carrying potential of smoke structures and symmetries by purposeful creation and real-time manipulation. Signal smokes were precisely sculpted with great skill, which also enabled the use of spatial and temporal frequencies and enumeration to encode information. 

\section{Results \label{results}}
The excerpts and quotes from 19th and early 20th century books, memoirs and newspapers collected in the Appendix comprise the primary results of this study. 
The key findings from this material are summarized in in tables  \ref{table1} and \ref{table2}. 
Some brief explanatory notes  ensue, and  these results are discussed further in section \ref{discussion}.  

\begin{table}[ht]\caption{\label{table1} Summary of historical and cultural information concerning smoke signalling, with  \# labels linked to corresponding excerpts  in  the Appendix. 
} 
\vspace*{2mm}
\small
\centerline{
\begin{tabular}{p{0.6\textwidth}p{0.3\textwidth}}
 \hline\hline
Description &  excerpt \# \\[1mm]
 \hline
Verified pre-colonial (pre-1788) sources  & \ref{j}, \ref{l}\\
&  \\
Spiral and fluid dynamics games played by children, training  &\ref{f}, \ref{m},  \ref{s}\\&\\
Puzzlement of white observers, their inability to understand smoke coding, mystery and secrecy of smoke signalling methods and codes&  \ref{a}, \ref{b}, \ref{dg}, \ref{k1}, \ref{l}, \ref{w}, \ref{z}, \ref{z1}, \ref{aa}, \ref{dd}, \ref{dd2} \\&\\
Tribal and locality smoke labels,  individual smoke signatures, named smoke signals & \ref{a}, \ref{k1}, \ref{l}, \ref{r}, \ref{s} \\&\\
Women smoke-signallers &   \ref{a}, \ref{b}, \ref{c}, \ref{o}, \ref{aa} \\&\\
Names of individual signallers  &   \ref{a}, \ref{c}, \ref{j}, \ref{l}, \ref{s}\\&\\
Relaying, repeating, replying &  \ref{b}, \ref{c}, \ref{d}, \ref{e}, \ref{j}, \ref{k2}, \ref{p}, \ref{q}, \ref{bb}\\&\\
\hline\hline
\end{tabular}
}
\end{table}

The linked excerpts  \ref{j} labelled in table  \ref{table1} tell us that smoke telegraphy was a well-developed technology by at least the 17th century CE. The book referred to is a factually-based account of the history of the Gunn-e-darr tribe of north-central NSW between c. 1650--1740 CE; specifically, the story of how Chief Red Kangaroo unified the widespread (over a district at least the size of Hungary) Gamilarray (or Kamilaroi) language tribes and clans. This history was recited by Indigenous Elder and historian of the tribe, Joe Bungaree, between c. 1890--1900 to Gunnedah police Senior Sergeant John Ewing and his son Stanley Ewing and recorded verbatim in real time by Stanley in hand-written notes\footnote{The original Ewing transcripts and documents are held by the Australian Institute of Aboriginal and Torres Strait Islander Studies, \href{https://aiatsis.gov.au/}{https://aiatsis.gov.au/}. Call Number: MS 1453.}. These precious transcripts eventually were written up by Ion Idriess in 1953 \citep{Idriess:1953}. In Joe Bungaree's words: `In every tribe there are men trained to remember. And so my grandfather trained me, and his grandfather before him.' 

A second pre-colonial oral history,  linked to excerpt  \ref{l}, is from 1770 CE. The book referred to is a carefully curated collection of Indigenous oral histories of Cook's voyage along the  east coast \citep{Rix:2024}. Of note is that each tribe had a `signature smoke', so that recipients of a smoke message would know which tribe had given its imprimatur to the content: `Each group would relay their own song and story within the smoke they sent'; in other words, an enumerated, ordered  group of smoke structures.   

Excerpts linked to \ref{f} and \ref{m} of table \ref{table1} indicate that mathematics was actively `lived' from childhood and embedded in Indigenous culture---not just an ordeal to be suffered only inside the school classroom. 
Learning about propellers and spiral symmetry and trajectories was a joyful and boisterous game played by children. When the time came for teenage initiates to begin training as smoke telegraphy apprentices their geometric, spatial  and number senses were already well-developed through active mathematical learning. 

Reported names of certain smoke signals,  excerpts linked to \ref{k1} and \ref{r} of table \ref{table1}, may be agreed abbreviations, perhaps across languages, of commonly used messages, analogous  to formerly widely known Morse code sequences, such as SOS and CQ (calling all stations) FB (good, fine business), SVP (please). A cry of `bobuckbee!' is more succinct and attention-grabbing than `I see a thin column of light grey smoke, someone is in distress'. 

Turning to table \ref{table2}, it seems certain  that the deployment of spiral streams was of crucial importance in coding. The explicit creation and use of enantiomorphic (right- and left-handed) spirals was described by Duncan-Kemp, excerpt~\ref{b}. The chirality of a spiral smoke stream is visually unambiguous and unfalsifiable (figure~\ref{figure2}). Taking advantage of symmetry is an  economical way of communicating a binary piece of information that must be more explicit than on-off: Duncan-Kemp describes how Mary Ann Coomindah makes and sends off a left-handed spiral, to signal that a white man was in the party. Since most people are right-handed, a left-handed signal would have been more difficult to produce, but  it was vitally important to communicate the presence of a white man clearly with no possibility of error. 
\begin{table}[ht]\caption{\label{table2}  Results summary of mathematical and engineered smoke signal qualities,  with  \# labels linked to corresponding excerpts  in  the Appendix. 
} 
\small
\centerline{
\begin{tabular}[t]{p{0.6\textwidth}p{0.33\textwidth}}
 \hline\hline
Smoke signal attribute or description &  Item \# in the Appendix\\[1mm]
 \hline\\
Spirals, curls, whirls, twists, and how  they are made  & \ref{a}, \ref{b}, \ref{c}, \ref{d1}, \ref{f}, \ref{g}, \ref{j}, \ref{k1}, \ref{m}, \ref{r}, \ref{s}, \ref{u}, \ref{v}, \ref{x}, \ref{z}, \ref{dd}, \ref{ee} \\
&\\
Left-handed and right-handed spiral signals  & \ref{b}, \ref{c} \\&\\
Discrete smoke puffs, balls, balloons, clouds, or  `breaks'  & \ref{b}, \ref{c}, \ref{g}, \ref{i}, \ref{k1}, \ref{l}, \ref{r}, \ref{s}, \ref{v}, \ref{bb}, \ref{dd}, \ref{dd2}, \ref{ee}\\&\\
Mushrooms, rings, and other smoke structures &    \ref{g}, \ref{i}, \ref{k1}, \ref{s}, \ref{t}, \ref{v}, \ref{y}, \ref{dd}, \ref{ee}\\&\\
Dictated complex messages and accurate readouts of smoke telegrams & \ref{a}, \ref{b}, \ref{c}, \ref{d1}, \ref{e}, \ref{h}, \ref{i}, \ref{j}, \ref{k}, \ref{k1}, \ref{l}, \ref{n}, \ref{o}, \ref{r}, \ref{s}, \ref{x}, \ref{z}, \ref{bb}, \ref{cc}, \ref{dd}, \ref{dd1}\\&\\
Numerically coded signals  & \ref{b}, \ref{d1}, \ref{e}, \ref{k1},  \ref{n},  \ref{o}, \ref{r}, \ref{t}, \ref{dd}, \ref{dd2}  \\&\\
Spatial and temporal frequency coding, ordering, position or direction   &\ref{b}, \ref{c}, \ref{n}, \ref{r}, \ref{t},\ref{v}, \ref{x}, \ref{bb}, \ref{dd} \\&\\
Distances over which signals are observed or observable  & \ref{b}, \ref{g}, \ref{h}, \ref{r}, \ref{z}, \ref{aa}, \ref{bb}, \ref{cc} \\&\\
\hline\hline
\end{tabular}
}
\end{table}

Discrete-analogue signals, or elements of a signal, consisted of timed and enumerated puffs or smoke-balls. Structures such as mushrooms and rings were used less often, and excerpt \ref{ee} mentions other, more exotic, structures being part of the code. 

In aggregate, the results  compiled in table \ref{table2} affirm that mathematical knowledge and practices were inherent in smoke telegraphy. 
\begin{itemize}
\item We recognise smoke signal waves as modulated analogue signals, where the carrier waves are the atmospheric flows and we have, by analogy with radiofrequency communication signals,  
\begin{itemize}
\item	amplitude modulation (AM), where the amplitude of the signal wave is tuned in response to atmospheric conditions;  
\item	frequency modulation (FM), where the frequency of the signal wave is varied, such  as timed smoke pulses and repetitive spatial patterns;
\item	phase modulation (PM), such as timed succession of a row of signals.
\end{itemize} 

\item Coding is used that can transcend language, using chiral symmetry and other structures such as rings, pulses with various frequency and length scales, enumeration, and colour. 
\item Canonically ordered operations are prescribed: message elements must be prepared and sent in correct order and  have their correct `place value'.  Order of operations, or \textit{bulkawilidjie} in old Mithaka dialect (excerpt \ref{c}), was fundamental to  every aspect of life in Aboriginal societies, including effective smoke-signal  communications.

\item The signallers had mastery of the technical, spatial, and numerical  skills needed to create and maintain high, precise smoke streams of specified geometry, colour, sequence, frequencies and count,    
and empirical knowledge of fluid dynamics. 

\end{itemize}

\section{Discussion \label{discussion}}

The following subsections provide additional explanations, informed speculations and historical background that evince a strong tradition of Indigenous mathematical knowledge. They also may inspire ideas  for development of  classroom and research activities. The references become (even more) eclectic, bringing  in physics, ecology, physiology, and ethics. This is not a bad thing: students from many disciplines must achieve at least first-year mathematics, and the Indigenous sciences always have been interdisciplinary.

\subsection{Systems, networks, and graphs}

With its heterogeneous physical and intellectual elements, smoke telegraphy is understood  as a line-of-sight, relayed communication system. Subsystems could spring up and operate on an as-needed basis (which by many accounts could be rather often during daylight hours). The signal hills, headlands and other high points of the landscape comprised the permanent infrastructure. 
Sporadic or mobile  local subsystems operated independently, too. 
As a network or graph, its nodes are not permanently connected like poles-and-wires systems. The strength of smoke telegraphy resided in the supporting knowledge.  

These  qualities may have sustained its effectiveness and longevity, and the story of smoke telegraphy could serve as a parable for modern systems and networks of many kinds that seem hyperconnected and more dysfunctional (or at least exasperating) than beneficial. There is some heavy-weight support for this take: A  mathematical result by \cite{May:1972} found that a network that is  too richly connected  (i.e., with persistent, high connectance) or  too strongly connected (i.e., with high average interaction strength of connections) has high probability to become unstable and collapse, and the effect is more dramatic the larger the number of nodes. 
May's result referenced ecosystems, and has been upheld since for engineered, chemical, physiological, economic, and social networks. 

That First Nations civilizations always had an excellent empirical understanding of May's probabilistic stability criteria for complex networks is borne out by their sustainable management of traditional fisheries \citep{Dargin:1976} and of landscapes using fire mosaics \citep{Sharples:2025}.  The systems properties of smoke telegraphy suggest that modelling of sparse,  fluctuating weakly or intermittently connected information-flow, driven-dissipative networks may give rise to interesting and genuinely novel phenomenology.

\subsection{Codes, encryption, and secrecy}
It seems likely that the coding elements were syllabic, or smoke `pictographs', like the characters carved or burnt on message sticks \citep{Idriess:1939}, and codes also used frequency and numerical ciphers. The use of a  `Morse code'  dot-dash (long and short smoke cloud) system is asserted or assumed by several sources (Appendix excerpts \ref{z} and \ref{bb}).  None of it was ever systematically deciphered by white observers. 
There always were secret codes known only between certain persons within and between tribes, too, and the 19th century frontier wars were a powerful incentive to further secret code development, to confound white police and military tactics  \citep{Kerkhove:2021}.  

In short, the details of codes were classified secret knowledge, which to this day has not been declassified by cultural authority, and there the matter should stand. Students may carry out some reverse-engineering based on information in this paper and although the outcome is unlikely to be `true', the exercise undoubtedly would engender profound respect for Indigenous know-how and  provide a memorable introduction to several core curriculum topics. 

\subsection{Symmetry, chirality and fluid dynamics}

From the results in section \ref{results} it is striking  that knowledge of chiral symmetry was `lived' mathematics of Indigenous societies. Physical manifestations of chirality (other than hands and feet) were well-known and part of everyday life. Whirlwinds and spiralling dust-devils (willy-willies) are seasonally common in many parts of Australia. `Birds of the whirlwind' are the vast flocks of budgerigars that rise from the trees in synchronized, spiral murmurations of green and gold \citep{Duncan-Kemp:1961}. Whirlpools are observed in wave, tidal, river or flood water races. The helices of land and sea mollusc shells are usually right-handed, but some species have left-handed helical shells. Climbing and tendrilled plants spiral around a support in a left-handed or right-handed sense. Propellers, which  are intrinsically chiral objects with blades set at a pitch or twist, were made by children from leaves for playing games (excerpt \ref{f} in table \ref{table1}). Some types of boomerangs were designed and  manufactured as chiral objects, similar to propellers,  for right- and left-handed hunters or sportspeople. By design, some sets of kinship rules were effectively non-superposable mirror image groups  \citep{Field:2025}. 
Language lexicons distinguish left-handedness from right-handedness, e.g., in Kamilaroi (Gamilaaray)  \citep{Ridley:1875} and in Guugu Yimithirr  \citep{Haviland:1998}. 

Chiral smoke streams were actively produced and exploited by Indigenous smoke signallers in Australia from at least the 17th century. In Europe, Louis Pasteur is credited with the discovery of molecular chirality in 1848, by examining the macroscopic morphology of the enantiomers of ammonium tartrate, which he separated with tweezers  \citep{Gal:2017}. The various types and subtleties of molecular chirality became important in the pharmaceutical industry after the 1950s thalidomide tragedies---only one of a pair of drug molecule enantiomers works as intended.  
In 1894 the terms chirality and chiral were introduced by physicist William Thomson, ``I call any geometrical figure, or group of points, chiral, and say that it has chirality if its image in a plane mirror, ideally realized, cannot be brought to coincide with itself'', 
and in the same year physicist Pierre Curie 
published a paper on symmetry in physics where he referred to handedness as enantiomorphous or nonsuperposable dissymmetry. 

It was not until well into the 20th century that mathematicians began to take an interest in chirality as a geometric property \citep{Weyl:1952}. In mathematics, the problem of chirality is not fully resolved to this day, for there is still no widely accepted, all-in-one definition of symmetry itself on which to base a general definition of chirality, although such a definition  has been proposed \cite{Petitjean:2007}. A pair of enantiomorphs is isometric, as well as being non-superposable mirror images, so we need coordinates to distinguish the objects  \citep{Gerlach:2009}. A modern mathematical definition of chirality (i.e., in spaces that are not necessarily Euclidean) is given by  \cite{Petitjean:2020}.  

Although chiral symmetry seems to have been known and expressed in crafts, pottery, sculpture and architecture for millennia (perhaps in pasta too), its use in smoke signalling may have been a first as a coding element in  a working communications technology. Chirality seems to be an instance of how empirical knowledge and effective technological usage often precedes theoretical or algebraic formulation.  

In the early 20th century, the vorticity of atmospheric flows (whether or not visualized by smoke aerosol particles) was recognized by \cite{Richardson:1922}, but that followed a long development through the late 19th century of the vector and tensor calculus that enabled the convenient definition and notation of vorticity as the curl of the velocity field, to a large extent by Heaviside---who was a telegraphist---and by Gibbs. The vector notation and formulation was controversial, and was opposed in the literature especially by Tait who was devoted to quaternions. 
However, Tait evidently decided that pragmatism was the better plan, for he carried out experiments---as mathematicians often did (Richardson, Fourier) and the smoke telegraphers would have done---with smoke vortex rings in which the rings interacted but retained their separate identity. (\href{https://mathshistory.st-andrews.ac.uk/}{mathshistory.st-andrews.ac.uk/}.) The production and interactions of smoke rings were well-known to Indigenous smoke-signallers, and were part of meaningful message coding (table \ref{table2}).  

\subsection{The frequency domain and the Fourier transform}

Of interest are the excerpts, labelled and linked in table \ref{table2}, that mention use of the frequency domain for coding. The number of repetitions of a smoke type or pattern within a unit of distance is its spatial frequency. This makes sense in smoke signalling because i)  the spatial frequency components of a smoke signal are meaningful for the coder and are easily set up, 
and ii) the Fourier transform of a visible wave signal  to the spatial frequency domain is a perceptual operation carried out by the visual system of the receiver. 

The Fourier transform 
takes a complex signal and breaks it down into its spatial or temporal frequency components. 
It is undoubtedly one of the most useful pieces of mathematics of all time. 
But the idea did not spring fully formed out of Fourier’s head in the early 19th century \citep{Fourier:1822}. There are \textit{always} antecedents. 
Histories of mathematics say those were the works of Euler, and Gauss, and James Gregory. 

Be that as it may,  when we probe non-Western histories we begin to suspect, then actually find, that many  `higher mathematics' concepts that usually are attributed as exclusive products of the development of European or Western mathematics over the last two or three centuries, and are taught only at university level, were already known empirically in diverse human societies, and expressed culturally. We should not be surprised. Two of these, chiral symmetry and space/time-frequency domain transforms, are presented in this article; others will be expounded in a forthcoming work. A sample concerning phase plane analysis (usually taught in second-year mathematics) may be viewed at \href{https://maths.anu.edu.au/news-events/mathematics-art/the-same-story}{https://maths.anu.edu.au/news-events/mathematics-art/the-same-story}, with accompanying text.

Canonical transformations are embedded in stories all over the world, most often between earthly and sky domains,  from the Seven Sisters to the Resurrection. They always served a purpose of encoding and transmitting complex knowledge and data. By Fourier’s time, with transformations conceptually ingrained in European culture  and increasing awareness of their power and potential among the educated elites of the Enlightenment, mathematics was codifying transformations of many kinds  in the language of algebra and complex analysis. 
Fourier was doing what he knew, working within his cultural space. The Indigenous smoke signallers were working within a cultural space too. The message must get through correctly, and the medium had a low signal-to-noise ratio, so it was a natural expression of human thought to construct the smoke waves in terms of meaningful frequencies.  

Those had to be picked out by the receiver. Vision is cultural as much as it is biophysics and physiology. What you see and interpret depends on your experience and training. High visual acuity was crucially important for the master smoke signallers, because of the complexity of the frequency patterns, smoke structures and colour combinations. 
For such elaborate signals to be effective, receivers must possess the visual acuity and cognitive training to rapidly perceive, parse and decipher them. 
 The eye's optical system forms a piece-wise spatial spectrum of the observed object, as well as its image on the retina 
 \citep{Campbell:1968,Braddick:1981,DeValois:1988,Westheimer:2001}. An experienced smoke telegrapher would develop a high visual sensitivity  to middle spatial frequencies and less sensitivity to noise at high frequencies and the relatively structureless, low amplitude, low frequencies. 
 
 Currently it is speculative, but not unsubstantiated from the above results and discussion,  that space/time to frequency domain transforms are universal to the human brain and their expressions in various media are universal in human cultures. We are forever seeking long-distance communications that we cannot access  through pheromones or other methods used by insects and animals. 

\subsection{Who were the smoke telegraphers?}

It seems that most older children and adults had an elementary knowledge of smoke signalling. Everyone would need to know how to send up and receive a danger or emergency signal within country. The professional master telegraphers, both women and men, would have had long training. They must  know how to code up a complex message given in words and send it off, perhaps after discussion between colleagues,  and receive and decode such a message accurately. They need to know how to keep a sharp edge to their smoke, elevate and direct their smoke, and read atmospheric currents. They must be able to make smokes of different colours and all the structures mentioned in section \ref{results}, including left-handed and right-handed spirals. No doubt there were local tricks of the trade too.   

To create and send and decode and relay such signals, over vast distances, the master smoke telegraphers had to be skilled practical mathematicians, coders, and fluid dynamicists. For that, there had to be a systematic base of intercultural education and training passed down the generations and across cultures with different languages.

From  the excerpts 
labelled  in table \ref{table1} we know the names and have partial biographies of several accomplished smoke telegraphers. Among these, we know most about Mary Ann (Maghroolara) Coomindah (c.1865--1929), from the memoirs of Alice \cite{Duncan-Kemp:1934,Duncan-Kemp:1952,Duncan-Kemp:1961,Duncan-Kemp:1968,Duncan-Kemp:2005}, who grew up  in the  early 20th century on  Mooraberrie station on Mithaka lands of the south-west Queensland river channel country. From these sources, Mary Ann Coomindah's remarkable mathematical and other talents were described in \cite{Ball:2025}. 

In excerpt \ref{j} of the Appendix, Burradella sends off two boys who have had \textit{some} training to Black Hill look-out, after making them memorise and repeat  a very complex message to give to the signallers to send out.  Burradella knows  the boys are not nearly skilled enough to actually smoke-up  the message!  

\subsection{Continuity and vitality of smoke knowledge} 
Smoke telegraphy was superseded by electrical cable and radio-wave telegraphy and telephony, but that did not happen overnight. 
In 1872 the overland telegraph from Port Darwin to Port Augusta was completed, connecting Australia to the world via a submarine cable from Java. That story is told as an exemplar cybernetic system in \cite{Bell:2025}, but the authors do not mention that an effective cybernetic communications system already existed. In fact news of the progress of the overland telegraph was smoke-signalled by the tribes along the route (excerpt \ref{r}). Smoke telegraphy persisted alongside the overland telegraph and across the continent for many decades. 

Electrical telegraphy was  rolled out within and between the cities from the 1850s, by the end of that century telephony was well-established in urban areas, and from 1905 radio-wave telegraphy was in use. Aboriginal and Torres Strait Islander communities were early and keen adopters of the new technologies where they were available, because of the traditional importance of long-distance communications and connectivity.\footnote{The party telephone line, where a single telephone number and line was shared by multiple subscribers, was still in use in remote areas until the 1970s. It was by no means unpopular, perhaps because it mimicked the broadcasting capability of smoke signalling. When the phone rang, or a smoke went up, everyone picked up and listened in! Everyone knew all the gossip! Conceptually, social media is nothing new: \cite{Duncan-Kemp:1968} could almost be describing smoke-tweets flashing `up and down the rivers'.} 

In far northern regions, beyond the reach of poles and wires and the broadcast range of radio transmitters, smoke signalling and telegraphy was still used, but knowledge of the finer points declined. It seems that the very real Japanese threat in WWII caused the military to suppress its regular use, presumably for national security reasons (see Appendix, excerpt~\ref{dd2}). After the war smoke signals became infantilized as a Scout-camp exercise that some readers' grandparents might remember, but as a sophisticated long-distance communications technology smoke telegraphy seemed to disappear from white public consciousness, from newspapers and books, and from the published white colonial histories.  

\newpage

But culture and knowledge are strong. In Aboriginal and Torres Strait Islander cultures smoke is still an important messenger medium. 
Today, the safe and expert deployment of smoke   is a vital element of cultural practices, and is used ceremonially as a sign of respect, in welcome to country ceremonies, and scientifically in traditional fire management and in preparation of foods and medicines \citep{Wettenhall:2022,Langton:2022}. 

To add further motivation to that stated in section \ref{section1}, the mathematics of smoke telegraphy also falls under the Diversity and Perspectives umbrella of Chiodo's and M\"uller's guide to ethical and responsible mathematics \citep{Chiodo:2025}. It encourages a diversity of people in mathematics through a diversity \textit{of} mathematics. Dismissing this knowledge as a mere footnote of history undervalues and alienates potential talent, limiting, down the track, the capacity of mathematics in productivity growth and innovation, which we know depends on diverse teams and perspectives. 

\section{Conclusions \label{conclusion}}

\begin{enumerate}
\item Australian First Nations smoke telegraphy applied mathematical knowledge  of chiral symmetry, frequency coding and decoding, and understanding of systems and networks.  
The mathematics was empirical, but not less effective for that ~\citep{Raju:2018}. There can be no stronger driver for accuracy than the imperative to send out a meaningful signal correctly, timely, and unambiguously, and to ensure this is systemic knowledge. 
\item This material  may be included in a cultural and historical introduction to first and second year tertiary mathematics and physical science courses where topics such as symmetries, groups, and the Fourier transform are studied. The systems properties of smoke telegraphy may inspire modelling and simulation problems to be  set up and analyzed with modern mathematical, computational and machine learning tools, and which may lead to novel and interesting insights in network theory. 
\item For Indigenous and minority students this study provides affirmation that mathematics is a valued part of their cultural heritage and identity and that Indigenous mathematical knowledge is rich and strong.  Students  are likely to find more meaning, value, and purpose in mathematics by seeing their own identities, cultures and histories reflected in the curriculum content and mathematical tasks given \citep{Leung:2020}, and in the cultural foundations; this is an educational \textit{right} that includes Indigenous students \citep{UnitedNations:2007}.  

\item A speculative, but not unreasonable, sequitur from this study is that development of mathematical ideation  in the human brain  may  occur at similar rates and levels in all societies, and that expression of those ideas is cultural in actuation but always serves a social purpose. Humans, it seems, are a social, mathematical species. 
\item Limitations of this study and further research: 
\begin{itemize}
\item This survey of published bibliographic and archival material is not comprehensive, and being confined to Australia  this is not a global catch-all study. Smoke signalling was practised in many parts of the world by many societies, but it appears to have been most highly developed by First Nations peoples in Australia.  This research may help inform Indigenous-led, ethical and culturally permitted reviews of unpublished material on smoke signalling. 
\item  Is there evidence that European mathematicians in the 19th century might have known something of Indigenous uses of smoke spirals? A fascinating question, posed by Prof. Falconer.  Where does this sit within the history of mathematics? Pre-colonial international connectivity had occurred for centuries between Arnhem Land and Indonesia and through the Torres Strait Islands, and post-colonial transmission of Indigenous knowledge to the UK is likely through letters and verbal transmission by returnees.  This is another project in history and archives, and there may well be dots to be joined. 
\end{itemize} 

\end{enumerate}

\bigskip
 
ORCID Rowena Ball \href{https://orcid.org/0000-0002-3551-3012}{https://orcid.org/0000-0002-3551-3012}.
The author reports there are no competing interests to declare.

\subsubsection*{Acknowledgements}

\textbf{Cultural Information:} This article contains Indigenous Cultural and Intellectual Property (ICIP) \citep{Janke:2021}, as it describes historical mathematical and scientific knowledge that is cultural and intellectual heritage of all First Nations of the Australian continent and islands and it references Indigenous histories. 
The collective ownership of and rights of ICIP holders to the knowledge described in this article are acknowledged and respected. The material in this article does not include culturally restricted knowledge or images or violate secrecy. 

I would like to acknowledge helpful and supportive conversations with Isobel Falconer and David Pritchard, of the History for Diversity in Mathematics Network, and with Ray Kerkhove at the University of Queensland. 

This work was supported by Australian Research Council grant IN230100053.

\renewcommand{\refname}{References\medskip

\small{\parbox{\linewidth}{\normalfont Many of the books listed have not been digitized (to date) and are not readily available outside Australia, except perhaps through online sellers of rare books. For that I apologize. Nevertheless, they are original source material and for that reason, and for their intrinsic interest, they are included. Websites referred to here and throughout the article were checked for accessibility on date of submission of revised manuscript.\bigskip}}} 

\bibliographystyle{plainnat}

{\raggedright \small
\bibliography{SmokeTelegraphy}
}

\small

\section*{Appendix} 

The following excerpts quoted from primary published sources are  
referenced in tables \ref{table1} and \ref{table2} of the main text by their numerical labels. 
Some edits have been made to omit or replace language that may be regarded as offensive, and to interpolate best-guess text where the print was so faded or damaged as to be unreadable in places.  Many of the historical newspaper excerpts are digitized at \href{https://trove.nla.gov.au/}{https://trove.nla.gov.au/}.

\subsection*{Published books, booklets, journals}

\setenumerate{leftmargin=0cm}
[Alice Duncan-Kemp's books, items \ref{a}--\ref{d}, are memoirs written from her diaries of the early 1900s] 
\begin{enumerate}

\item \label{a}  \cite{Duncan-Kemp:1934} 

p. 198: Usually, the exact meaning of smoke signals, whether dense clouds, pale spirals, columns, or breaks, depends on circumstances. A pale wisp might mean a warning regarding a trespass of game laws; two columns, victory to a warring party. If a party were out hunting the same signal would mean that game had been found. 

It was an education to watch the women gather dry leaves or grass and place them so that one pile would send up a faint buff coloured smoke signal, while another pile would emit a dense cloud. Sometimes only a faint wisp escaped, barely discernible, but the ever-watching blacks would see and interpret it immediately, leaving the puzzled white man to wonder how the blacks knew so much about the stockman’s movements. The code for smoke signalling is intricate and very elaborate. At a glance a black can tell whether a white man is mentioned in the streak fluttering faintly on the horizon, and who and where the man is. One day Mary Ann shaded her eyes with her hand and said: ``That one Andy Clinton, Monkira boss, leg cut with an axe.'' And it was true although a week passed before any word came through. 

A signal known and dreaded by all tribes is the \textit{ko-to-unjie}, or death signal. It had a tail or wing to it; the more important the person, the bigger and longer the smoke and wing. 

p. 203: A black smoke cloud twisted and twirled, now heavily in pulses, now lightly in pearly columns. For half an hour we watched the smoke ``talk'' from the \textit{Kibulyo} camp, sending out a  code message to be relayed by surrounding tribes. Suddenly a deep black-winged column of smoke shot up beside a pale grey streak. It was the \textit{Koto-unje}, a death signal.

p. 205: The interpretation of these signals often puzzled us. It left us convinced that the blacks’ system of bush telephony has never been, and never will be, fully revealed to the whites.

p. 207: A tell-tale wisp of smoke floated over the western line of sandhills. Mary Ann and Sally were galvanized into activity. It was a telegram to say that Tommy, the Monkira horseboy, would be at the station camp that day. About four o’clock Tommy and Pituri Bob arrived.

p. 279: Small smoke clouds rose from the burning totems.  The \textit{kobi} [medicine man] stared long at the curling wisps, seeking a trace of a face or some sign that would warn his tribe of those responsible for such evil. Slowly a mighty hunter took shape under his gaze---Goolburra! 

\item \label{b}  \cite{Duncan-Kemp:1952} 

p. 11: Excitement ran high. The women scanned the western horizon for smoke signals, the bush telegraph which relays in some intricate secret code tribal news to the outside world. 

p. 19: With a pile of green leaves sizzling in long clouds of brown-grey smoke, Mary Ann and Bogie tapped out a message to the Waker-Di camp 12 miles distant. Now short, now long, and short again, came the smoke stream. With a light bough they regulated the waves until they became thin columns and smoke breaks, or delicate spirals with a peculiar left to right curl which is said to mean `white man', instead of the more conventional right to left pointed spiral.

p. 22: Suddenly a great hush fell upon the party; thirty pairs of eyes were glued on the southern horizon. A black smudge wavered over the sandhill’s crest, then up shot a thin dark streak into the air. 
Up flared the signal, now light, now dark the spirals wavered in the breeze, then in short dense breaks alternating with dark columns the smoke waves came to us. That an important message was being transmitted I could tell by the rapid change in expressions. 

p. 175: It is strange how the blacks can tell the difference between a white man's smoke and that made by one of themselves. Shortly after the white man's signal another smoke rose over the sandhills: ``Old man belonging to Coonundhurra tribe come up tomorrow,'' said Bogie. 

p 185: From the sandhill peak we could see three smoke signals, one grey-black puff ball rising between two pearly-grey streaks---the Eaglekawks' signal used for telling surrounding tribes that a chief was on a walkabout to their camp. 

p. 192: About five o'clock that afternoon a pale smoke streak appeared above the north-west timber line---a telegram from the Pelican people to advise ``all those interested'' that Dingyerri with her son had arrived safely. 

\item  \label{b1}  \cite{Duncan-Kemp:1961} 

p. 137: The Murranuddas at Monkira camp read the news of Big Peter's death in the smoke-signals that wafted across the sandhills. 

\item  \label{c}  \cite{Duncan-Kemp:1968} 

p. 160--161:``Even smoke signals have their own \textit{bulkawilidjie} [protocol, place value, order of operations],'' explained Kimurni, ``Look!'' Wreaths of smoke coiled themselves into hoops and twisted in spirals as messengers transmitted their talk in smoke telegraphy from tribe to tribe up and down the rivers. It was all part of a vast network of an Aboriginal communication system.
Other tribes relayed to still farther back tribes the news of the moment. 

A deep black V-shaped column sprang up as the grey spirals died down. ``Old Murran of the Poorakees has put the \textit{bulkawilidjie} code into action,'' whispered Kimurni, "But the \textit{koto-too-eera} [smoke signals] get through just the same. Now the Wahlitchas will have to pay the Poorakees a big fine." It was noticeable that smoke signals were off the air for days afterwards. 

\item \label{d}  \cite{Duncan-Kemp:2005} 

p. 227: There, fluttering above the skyline just clear of the sandhills was a ragged winged smoke. It rose high, black and menacing over the countryside, trailing its wing like an injured bird: \textit{Koto-unjie}, the messenger of death. 

p. 261: One mid morning the southwestern horizon was alive with smoke signals, thin blue-grey wisps that merged into a continuous network over the sandhill’s backbone.

\item \label{dg}  \cite{Gunn:1908}

p. 200: Two tell-tale telegraphic columns of smoke, worked on some Aboriginal dot-dash system, had risen above the timber, and their messages that told them we were coming had also been duly noted down at the Red Lilies and elsewhere, and acted upon. 

p. 207: Dan regretted we had no way of letting Jack know. We were not aware of the efficiency of smoke signalling. Next day we rode across the salt-bush plain, to be greeted by smug, smiling old blackfellows and Jack who, to Dan's amazement, already knew. Dan believed they must have supernatural powers, until, at the homestead, Nellie revealed the smoke-signalling secret.  

\item \label{d1}  \cite{Hassell:1975}
[An account of the author's interactions with the Wheelman tribe, near Jarramungup, about 80 km inland from Bremer Bay, in the 1890s.]

p. 122: On an excursion to the gorge,  Greton suddenly saw two thin spiral coils of white smoke curling up to the sky, then further on two more. She said it was a signal to go back.

p. 193: After travelling three days the party saw several fires, so decided to make fires and signal to them. Towards evening their signals were replied to and the following day both parties travelled in a direction so they should meet, making smokes as they went to show their route.

\item \label{e} \cite{Howitt:1904}

p. 720: Ngarigo signallers would go up onto a hill and fill a rolled sheet of bark with dry grass; by setting fire to this a column of smoke would ascend and be responded to in the same way. 
Kamilaroi signallers would select a tree with two or three pipes in it. Smoke issuing from two pipes might mean peace and three war. If the tree had too many holes, they would be plugged. 

\item \label{f}  \cite{Idriess:1933}

p. 28: In north Qld the Aboriginal youngsters would make their own fun. After making a dense black column of smoke, boys and girls would  `spear' long leaves into it. Immediately the leaf was caught within the rising current of hot air it would spiral violently upwards; the winner was whoever could make his or her leaf rise highest. When a score of leaves were spiralling up in a race the excitement was great as all eyes gazed upwards and all voices cheered the leaves on. 

I first saw my `aeroplane propeller' long before I saw an aeroplane, when those children would shape two long leaves, then centre them using a thorn for a pin. When  `slewed' into a smoke column it would right itself in the centre of the hot air-current and whirl itself swiftly upwards, just as propeller blades whirl round. Surprising things the youngsters could do with leaves of varying textures, length, breadth, and shape, and with or without a heated air-current, for they could sometimes find a natural current that seemed to an ignorant white man to be invisible. 

\item \label{g} \cite{Idriess:1935}

p. 176: The Aborigines signal in big or small puffs, long thin stems, mushrooms, and spirals. They can send smokes to different heights, and in varying colour, speed, and density. They can select air currents rising from the plain to lift it high  above the lookouts.

Chugulla and Nemarluk's men made the smoke signals in this way: a  heap of dried leaves that give ``long distance'' smoke were piled cone-shape and lit. Then a mat of weaved grass was thrown over it. Under the edges of the mat, just sufficient air is allowed in to keep the leaves burning and forming the greatest volume of smoke which slowly rolls over itself under the mat. At a grunt from the leader, the mat is sharply withdrawn with a peculiar whirling motion that draws air under the mat and helps to shoot up a puff of smoke. The mat is immediately replaced and another puff formed. A still day with clear air is best for signalling. Then the air currents seem to rise up from the plain to lift it high into the air above the look outs. The visibility too, is such that a column, a balloon, a cloud, a mushroom, a spiral can be seen at great distances.

\item \label{h}  \cite{Idriess:1939} 

p. 24: From a peak across a valley behind us rose a snaky wisp of pale blue smoke. As it climbed into the still air it swelled stronger into a black column almost stationary.
``Bad signalling country,'' nodded Laurie. ``Too many ranges limit the vision. All the same, in a few hours the Aborigines for a hundred miles around will know that the patrol is at Mount Hart.''

p. 45: ``You can't beat Aboriginal smoke-talks. Only wireless could beat `em, but we've got no wireless out here. I often hear news away out in the hills that Derby doesn't hear for a month.''

\item \label{i}
 \cite{Idriess:1941}  

p. 87: From the look out rose a puff of black smoke followed by another then another, then a ring of smoke with a black puff shooting up straight through it. Slowly the smoke drifted away. They halted, staring at one another. This was an urgent danger signal.

\item \label{j} \cite{Idriess:1953}

p. 61: Several men ran for the nearest hill to summon the scattered hunters by a smoke signal. Faintly, from away over towards the Whispering Belahs, came floating an answer from Jerrabri's men. They in turn would signal Kuribri's men at Tambar Springs, while those fishing down the Namoi would smoke-signal the look-outs on distant Ydire and on Booroobil Rock. 

p. 168: Burradella acted swiftly, and definitely.
``Bunnadunne!'' he said. ``Race to the Black Hill. Tell the look-out to smoke-signal Carrowreer and the Red Chief away on Breeza Plains: {\sc Enemy been heard talking on river crossing between Burrell Lagoon and Mullibah  camp. Women and children being rushed to the Secret Camp under guard}.'' And the swiftest runner of the young warriors was speeding away towards the look-outs. 

pp 174--175: ``Kuppa! Burrai!'' Burradella said sharply to two of the older boys. ``You both have been trained in smoke and fire signalling. Now, run to the Black Hill look-out, tell the three men there to send smoke signals to the Red Chief saying, {\sc Cassilis enemies. Four seen on river. Looks like war party camped close.} Understand?''
The boys nodded. 
``When they have signalled, tell them to signal the two look-out men away on Carrowreer to repeat that signal. After which the Carrowreer men are to make their way to me here. Understand?''
The boys nodded again. 
``Right. 
You two stay on look-out duty in case the Red Chief signals. If so, one of you hurry with the message to me. Now run!''

p. 181: From the Porcupine Ridge look-out a long finger of black smoke was shooting up into the clear blue sky. As they watched, the smoke column changed to white. After an interval, it changed to blue.  It was duplicated swiftly from distant Carrowreer. Then a big ball of black smoke shot up, followed by a small ball. Then, at intervals, arose black, blue, and white smoke curls. There was an interval, then the black column shot up again. The signal was being repeated. ``Kuribri'', said the Red Chief, ``tell me what that smoke talks to you.'' ``{\sc Enemies on Namoi crossing. Come quickly. Secret camp. Burradella.}'' answered Kuribri. The Red Chief waited until the signal was being repeated then said, ``Tell me what the smoke talks to you, Yuluma.'' And Yuluma read slowly the same message. ``And it talks so to me too'', said the Red Chief. ``I shall stay here long enough to signal Boobuk and his men.'' And he kicked one of the fires into a blaze and hurried to pull armfuls of grass and reeds, dry and green, for the different smokes they would give out. 
Soon his smoke signals were rising up.

\item \label{k} \cite{Kral:2020}

p. 32: Generally, on leaving and entering waterholes during \textit{yalatja} [short hunting trips], families would light small fires [for signal smokes]. By lighting a fire they were sending a message: 'Yes, come, water or game is in abundance here' or 'Yes, I am at the Warakurna soak and your sister is here with my aunty and uncle'. 
 
 p. 33: In the early days everyone travelling around, all the time looking for rain and smoke. 
 
 \item \label{k1} \cite{Magarey:1894}
 
There are several forms of signalling smoke, and amongst others the following are in regular use: (a) A pale-hued slender column, (b) a pale-hued heavy column, (c) a slender column of black  smoke, (d) a heavy column of black  smoke, (e) a spiral-coil form of pale and dark (each) smoke, (f) interrupted or intermittent smokes, e.g., in puffs, balloons, balls, parallels, etc, (g) groups of smokes, one, two, or three of any one form, or of several forms at times merged. 
 
There are combinations of form from the same smoke-column by means of sudden wavings of sheets of bark, causing side puffs of the smoke to rise on alternate sides of the producing blaze. By skilful manipulation of the smoke-pillar two columns rise in the air parallel to each other, one column pale-hued, one column dark. 
 
 The spiral form is obtained by using bushes or clumps of grass that have been twisted into a denser mass  and lighting the top. The fire rushes with intense velocity through the growth, and the smoke rises in swirls and coils into the upper air. The same effect is produced by the manipulation of a large skin or  rug; two natives standing on opposite sides of the fire and giving the skin a circular movement on an inclined plane, the smoke is formed into coils as it ascends. 
 Balloons or balls of smoke are  produced  to secure a succession of five or seven at regular distances in a column. 
 
 Festoons of smoke are employed upon occasion, usually used by a native moving rapidly who wishes to convey a quickly passing message to his own people. A string of grass is made and hung upon the boughs of a convenient tree, then lighted, and the signaller speeds on. 
 
 In Victoria, on the Darling, and in the Northern Territory, hollow trees are used for signal smokes. Green or damp material is thrust into the hollow, dry fuel is placed beneath, and the hot blast rushing through carries the dark smoke up through the tube and into the air in the form of a thin, regular, but far ascending dark column of smoke. 
 
Owing to the disinclination of the natives to yield up information in their possession concerning smoke signals, it is difficult to learn much as to the meanings attached to the various forms. 

The meanings of spiral coils are curious.  For Powell's Creek Tribe, spiral coils of thin pale or dark smoke ``Mullagar Winlabardim'' mean ``All about, come quick, plenty of kangaroo.'' Similar coils of dense dark smoke ``Umbarunnie'' mean ``Two men come quick, help carry game.

 Whilst the survey of Port Darwin country (1869) was in progress, the officer in charge was one day informed that war signals were being made by the natives. Upon ascending the hill to investigate, two spiral coils of light smoke were observed. Skins held by two natives were kept turning with a circular motion in an inclined plane over the rising smoke so as to cut the column at each revolution of the skin, and to give a spiral form and motion to the smoke as it rose. 
 
 \item \label{k2} \cite{Mitchell:1839}

By the smoke which arose from various parts we perceived that the Aborigines were watching our proceedings.  
When we reached the head of the highest slope,  a dense column of smoke ascended from Mount Frazer, and subsequently other smokes arose, extending in telegraphic line far to the south, along the base of the mountains; and thus communicating to the Aborigines who might be upon our route homewards the tidings of our return. These signals were distinctly seen by Mr. White at the camp, as well as by us. 

\item  \label{l} \cite{Rix:2024} 

pp. 18--19: Gunnai man Wayne Thorpe  tells that it was his ancestor, Boondjil Noorrook, who first saw the \textit{Endeavour} off Point Hicks in April 1770. His Elders said to him: ``Send a smoke signal and a song and forewarn all the people up along the coast.'' Boondjil Noorrook sent off the smoke signals to communicate with all the different mobs along the coast.

 Other peoples along the coast, having sighted the \textit{Endeavour}, also relayed smoke signal messages, forewarning people that there could be danger.  Each group would relay their own song and story within the smoke they sent. In Joseph Banks' and in Lieutenant Cook's journal they both write about  where they sighted these smoke signals and fires. 
 
The people along the coast knew the  different types of smoke, such as ceremonial, cooking fires, cultural burning and signal fires. The fires that accompanied Cook all along the coast were carefully managed signal fires situated on headlands. The journals of the voyage mention smoke consistently along the coast, but the British didn't understand that what they were seeing was a warning signal, alerting others to their presence.

p. 28: Worimi man Lawrence Joseph Perry  explained how First Peoples could communicate with each other over great distances without leaving their own Country, through smoke signals. He said that a fire would be lit on a specific headland or high place, and green branches thrown onto it to create a cloud of blackish-grey smoke. He also said that the heavier branches could be used to alternately cover the fire and release smoke puffs into the sky.  

pp. 294--295: Waubin Richard Aken tells that before Cook arrived in the Torres Strait there had been fire and smoke signals from Adhai (Crab Island, on the western side of Cape York), letting all the clan groups  know that a strange ship was travelling along the east coast. He says that smoke signals and the message stick were ``a blackfella internet''.

\item \label{m}  \cite{Roth:1897} p. 130:  

Games of  ``Smoke-Spirals'' are played in several north Qld tribes.  Any leaf, or even a mussel shell, by means of a peculiar motion of the wrist and arm, can be thrown  into the smoke rising from a fire so as to ascend with it like a spiral. Different ways of holding the article are shown in figure \ref{figure3}: (a) the wrist is rotated outwards as the fore-arm is jerked sharply forwards and downwards the object leaving the hand just at the completion of the movement;  (b) Another method is to hold the leaf between the fourth and fifth digits, and with a motion of the fore-arm similar to that in (a) rotate the wrist inwards; the Kalkadoon speak of this game as \textit{pi-ri-jo-rong-o}; (c) An easier method  is to make the  extended fore-finger of one hand act as a sort of spring on the other, the leaf at the moment of release being shot at an angle into the smoke.
 
 \begin{figure}[ht]
\centerline{\includegraphics[scale=1]{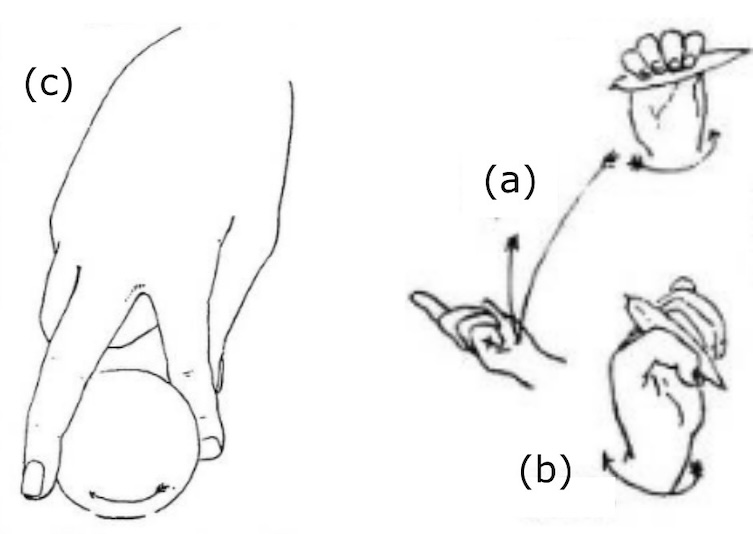}}
\caption{\label{figure3} \footnotesize Ways of holding objects to be rotated into a smoke stream so as to initiate a spiral trajectory, figure adapted from \cite{Roth:1897}.} 
\end{figure}

\item \label{n} 
 \cite{Ruhen:1948}. [A work of fiction, based on and embedding  true stories, that begins in 1902.]

When they had ridden on perhaps for two miles Edrington looked back. A thin column of smoke was rising to the zenith. 
``Is that a signal fire, then?'' Marriner asked, and his companion nodded. 
``He's telling his friends that he has met us but come to no harm. What they can read into it is amazing. His friends know he is out here. They can tell there's a new element come in. He's put up  smoke behind us, so they know he's met us. He hasn't changed his plans or he'd put up two smokes or three in a line in the direction he is now travelling. If they rose a long way apart in a short time they'd know he was travelling fast. There's only one smoke, so he's still doing what he was before and the inference is that he is untroubled by our presence, that we are friendly. There's another inference too, that he found us in no need of help from him. You'd be amazed at what you can read into that smoke.''

\end{enumerate} 

\subsubsection*{Newspaper reports, arranged in dated order}

\begin{enumerate}[resume]

\item \label{o} 
The Courier (Hobart, Tas.) 14 Apr 1854  p. 2.  
The Female Sealers

The Aboriginal women ascended to the top of a small hill and made smokes as signals to the natives on the mainland that they had taken some 12 seals. The smokes were soon answered by smokes on the beach.

\item \label{p} 
The Australasian (Melbourne)  28 Mar 1885 p. 43.  
Bush Life Forty Years Ago. By Myall. 
 
The Aborigines divided into two parties, and during their retreat communicated by smoke signals, which is done as follows: A fire is made and damp bark is piled on in such a manner as to leave a funnel-shaped
opening at the top; the dampness of the bark causes a white smoke or vapour to ascend in a thin column; the meaning of these signals is of course pre-arranged.

\item \label{q} 
Independent (Footscray, Vic.)  29 May 1886 p. 3. 
Bush Travelling in South Australia.

In 1862, my Government appointment  necessitated my travelling a route which had a radius of 200 miles from Port Augusta. I observed the native trade fair in which ochre and  bush tobacco (pituri) are traded. The South tribe know the time to a day of the arrival of the North tribe by smoke signals from hill to hill.

\item \label{r} 
Evening Journal (Adelaide, SA)  28 Oct 1893 p. 6. 
Smoke Signals of the Aborigines.

The triumph of the Australian Aborigines'   intelligence is exhibited in the realm of their telegraphic smoke signalling. Here they are without rival. 
Promptly, accurately, intelligently, the signaller hangs the cable message in the sky, and knows that no break in the line will stay it from its distant goal. It is very certain that each member of a tribe has
his own code-smoke by which may be conveyed any intelligence concerning him.
Each locality of interest in the tribal territory also has its own distinctive signal.

From Captain Cook we obtain the first reference to these Australian signals. On April 20, 1770, the Endeavour, having rounded Cape Howe, was steering northward, and Cook for the first time sees the warning smokes of the tireless, watchful natives. Later, Flinders, Eyre, Mitchell, Sturt, and Stuart were greeted with these signals in all their travels around and across the continent. 

In the Powell's Greek tribe, a thin pale-hued column of smoke  named ``\textit{bobuckbee}''  is a distress call. On sighting a kangaroo a hunter will fire material to hand in rapid succession as he pursues the game. By this means the direction and course taken is indicated, and others may join him. A slender column of dark black smoke   named ``\textit{Hugullo}''  indicates that a messenger is on the way. They can  read  smoke signals from from Newcastle Waters, fifty-six miles distant. 

A large dense column of dark smoke  is unquestionably the long-distance signal of the Australian continent, and is seen at about sixty miles distance, on account of  the moisture conveyed by the great uprush of smoke forming a cloud at the top of the column. This was used by the  Barrow's Creek tracker (from whom the information was obtained) after travelling eighty miles. His signal was observed and answered by his mate from his starting point in twenty minutes. 
In MacDonnell Range Tribe and Tennant's Creek Tribe, the thin black smoke column is named ``\textit{Quoorta}'' and means ``coming back'', used to signal  an unsuccessful search for water. 

Ball or balloon signals are produced by collecting the 
smoke in a skin inverted over the rising smoke. When the bag is full one of the assistants releases  it, allowing the smoke to escape in the form of a dark ball. This manoeuvre can be repeated with great rapidity and regularity.

A spiral coil form of smoke, which is in frequent use, is produced by firing a twisted bush, or manipulating the column with a skin, or twirling a bush round and round the rising smoke. 
Intermittent smokes are produced by interrupting the rising column, alternately rapidly and at long intervals. Sometimes two or three or five or six smokes are raised a few yards apart and manipulated as the nature of the message requires.

The Port Darwin peoples make a complex signal by the use of a
large sheet of bark stood on end near the smoke column, and by a quick downward movement the sheet is flapped towards the smoke, driving it out of
the perpendicular and to one side. By such manoeuvres on alternate sides, puffs conveying a
special code-meaning are repeatedly produced.

During  the Overland Telegraph Iine construction, 1870--1871, a camp was formed at Leichardt's Bar, Roper River. The Aborigines of the coast kept those at camp informed of the presence and arrival of vessels off the mouth of the river by signal smokes, which were relayed along the river up to the camp, a distance
of sixty-five miles. 

I have seen several column signals start up at great distances apart and 
approach each other by stages, eventually meeting  as  two
groups travel from opposite directions, heralding their approach and direction. 
A husband and wife  will separate for hunting early in the day, and by special form of smoke signal will indicate intention to rendezvous at a specified point to camp. 

When two workers were killed in the railway accident of
Monday, July 21, 1890, the news  reached Charlotte Waters Telegraph Station
the following day. The Aborigines there conveyed the facts promptly and accurately by smoke signals to Crown Point Station, fifty-six miles distant. The news was smoke-signalled on, and was known at Johannesburg Mission Station, 26 miles north-west, very early on Wednesday. 

Mr W. P. Auld, a member of Stuart's expedition 1861--1862, relates that on the
return journey (from the Indian Ocean) the party reached Mount Margaret where they met a station-hand who had learned six weeks previously that the expedition was returning. The number of horses and men was also known, the news having been telegraphed 600 miles by Aborigines' signal smokes. 

The signals vary not only in form and hue, but their significance also may depend upon
the hour of the day, frequency and locality. There is ample variation of signal to enable difficult and complicated messages
to be sent on by means of these smokes to friend or foe, far or near. The whole mystery of the code-meanings of these Aboriginal smoke signals yet remains unsolved.

\item \label{s} 
The Express and Telegraph (Adelaide, SA)  24 Nov 1893 p. 3.   
Report on meeting of The Royal Geographical Society. (See also \ref{k1}.)

Mr. A. T. Magarey gave an address on ``Smoke Signals of Australian Aborigines.'' Those trained in smoke signals could play on a fire just as a man would manipulate a keyed musical instrument. There is a code meaning attached to the various characters assumed by the smoke. Every tribe had its own particular code, while there was also a long-distance inter-tribal code. He gave instances of spiral smoke signals. 

Mr. Mark Wilson, an Aboriginal native from Point Macleay, described the smoke 
signals of the lakes and Murray River tribes. If a corroboree was to be held, a fire sending up a dark column
of smoke would be lighted at Mannum as a signal to the Point Macleay tribe, who would signal
the other lake tribes and all would assemble at the Wellington bora ground. He described the signals for war, death and losses in the bush, 
how fires wore built for night signals, and how puff smoke balls were made by holding the smoke under a skin or board then suddenly releasing it. 

Surveyor-General Mr. Goyder mentioned that when he had been camped with a party near Port Darwin, he saw spiral smokes produced by turning skins around, straight and small ball smoke signals, and flashlights made from the inner bark of certain trees.

\item \label{t} 
 The Australasian (Melbourne, Vic)  11 Mar 1911 p. 55. 
Smoke Signalling. By E. S. Sorenson.

Bright shells and glass-like substances were used by the
Aborigines for transmitting messages over great distances, however 
the most common method of signalling
was with smoke, the significance of the position, density, and number of the ascending columns being widely understood.

\item \label{u} The Queenslander (Brisbane, Qld)  14 Oct 1916 p. 41.   
Aboriginal Reminiscences by W. Clark.

The Western tribes used smoke signalling. Often when out on open downs country on the Comet
and Nogoa I have seen a spiral column of white smoke suddenly shoot up to
a great height on some distant range.

\item \label{v} 
Smith's Weekly (Sydney, NSW)  Sat 7 May 1921   p. 19. 
Smoke Signalling
 
The Aborigines' method of telegraphing over long distances was with smoke, which
they manipulated with wonderful skill, the position, density, colour, straight or spiral form,
and number  being part of the code. The atmosphere was kept as clear of
cooking-smoke as possible, so as not to interfere with the operations of the telegraph department.
The spiral form was obtained by moving the cloaking material around the rising smoke (figure \ref{figure4}). With accurate use of materials, numbered black puffs could be sent aloft alongside a column of white smoke. A hollow tree-trunk was sometimes used, the column
being broken as often as desired at the base. 

\begin{figure}[ht]\centerline{\includegraphics[scale=1]{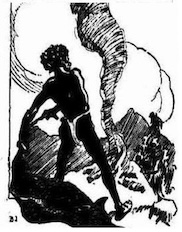}}
\caption{\label{figure4} \footnotesize Raising a spiral smoke, drawing from \ref{v}.}
\end{figure}

\item \label{w} The Daily News (Perth, WA)  18 May 1921 p. 3. 
Smoke Signals. 
 
While on the hills I saw a smoke signal on the Burtville hills that was remarkable. The column of smoke was divided into two branches,  but how the Aborigines from the township could tell it was Billy
beat me altogether. Their system of bush telegraphy is guarded by strict secrecy, and their efficacy in deciphering bush cryptograms is almost uncanny.

\item \label{x} 
The Queenslander (Brisbane, Qld)  10 Sep 1921 p. 11. 
  Aboriginal Smoke Signals. By Eli Rhys.
  
From observation of these signals Cook ascertained that Australia was inhabited. Flinders in 1823, then exploring Spencer's Gulf, inferred from the numerous columns of rising smoke that the country was well
inhabited and plentifully supplied with food. 

When need arose the Aborigines sent into the air messages of welcome or warning, invitation or defiance, mourning or rejoicing. Friends were invited to fish or hunt. Warnings were sent that a well had dried up. The presence of white strangers  was spread with amazing speed. 

The signals employed  may be continuous or intermittent, may rise as a single column or take a spiral form, they may give forth white or dark smoke. 
A slender column of pale smoke, in the central part of the Northern Territory, is a sign of distress. Farther south it is the warning sign that strangers or explorers are about. At Powell's Creek a heavy column of white smoke  means that a friendly tribe is coming
to talk. In other parts of central Australia it is a long distance warning of a death. A slender column of dark smoke is usually an invitation to a conference. By firing a large quantity of fuel, overlaid with green bushes, a column of dense smoke is formed rising to 1500 or 2000 feet. This is often a signal that strangers have crossed to the neighbours' territory. It interferes with atmospheric conditions and is often followed by rain. 

\item \label{y} 
The Herald (Melbourne)  3 Jun 1922 p. 11.  
Warrigal Wireless. An Aboriginal News Service.

These boys could  signal to each other from long distances with smokes. One would strip a ring of bark
four feet long off a tree, light a wood fire, throw some green leaves on it, and by holding the bark like a chimney over the smoke, send a column of smoke straight up. By manipulating the bark he distorted the smoke into all
sorts of different shapes, which could be read quite intelligibly by the natives to whom he was signalling. 

\item \label{z} The Sun (Sydney, NSW)  4 Oct 1925 p. 17.  

Long fingers and spirals of smoke stretch suddenly into the air, die away, and are followed  by others. To the Aborigines they are
pregnant with meaning. They constitute the ``bush wireless'',  by which they communicate over great distances.
  
From the bottom of the valley came a long spiral of smoke which was succeeded by short and long bursts of smoke.
Some days later O'Shaughnessy arrived at Newry Station. He was greeted by
the manager, who astounded him by saying that the Aborigines at the station had learned two days ago 
by bush wireless that he was on the road. They knew that he was wearing
khaki trousers, that he had a pack on his back, a gun in his hand, and that he carried a water bottle and was
accompanied by a dog. The message had been sent 50 miles from Aveurgne Pass and he 
heard later that  messages are received up to 80 miles by this method. 

In order to obtain the spiral effect, green
wood is used for the fire. A continuous stream of smoke is first sent up in order to warn the distant people
that a message will be sent. Then a blanket is used to cover the smoke at intervals, to obtain long and short bursts of smoke, which float into the air carrying a meaning. 
Particulars of the Aborigines' ``Morse code'' O'Shaughnessy  could not obtain, but he found out that each symbol means
a word or phrase. He also discovered that, like the Royal Navy, the Aborigines change their code
frequently. 

\item \label{z1} The Mercury, 12th Jan 1926 p. 5. Michael Terry, Bililuna Station. Smoke Telegraphy. 

The Aborigines seem to have some mysterious code whereby columns of smoke tell stories whose lengths and detail is only believable when actual instances come to note.

\item \label{aa} 
Telegraph (Brisbane, Qld), 27 July 1926 p. 20. 
Aboriginal Wireless. 

In the early days the pioneers were greatly puzzled to know how it was that,
wherever they went, their arrival was anticipated by the Aborigines. The riddle remained unsolved until the
white man realised that the Aborigines had a wonderful system of smoke signals by which the doings of the bush were broadcast for hundreds of miles.
The fact that these smoke were well understood by all is proof that practically throughout the continent, a common code prevailed. 

It is the women who collect the material and carry out  all the real work connected with it. Sometimes the men do the signalling, but they usually leave it to the women, who are equally expert. One can imagine that, when the first fleet was sighted in Sydney Harbour, the women of the great Gamaragal clan sent up smoke signals warning of the event, which, within a few hours, would have been repeated hundreds of miles.  

\item \label{bb} 
Toowoomba Chronicle and Darling Downs Gazette (Qld)  31 July 1926 p. 3. 

The Aborigines have evolved a thoroughly efficient system of conveying messages over vast distances by means of smoke signals, which are as clear in their meaning to them as an electrical telegraphic message is to an operator.  It is a version of Morse code, and must have been in common use for thousands of years. When an important smoke signal went up, it was repeated by other tribes, and so it would travel on till it
reached districts hundreds of miles away.  

The signaller employed different methods  for different occasions. 
A method which required a good deal of skill and ingenuity,  was to select a hollow standing tree with a hole at the base. A fire would be made in the hole and green leaves thrown on it, when the draught created by the fire would send the smoke racing up through the hollow trunk to a great height. When it was desired to break the smoke column a sheet of bark or a kangaroo 
skin would be thrown over the fire and rapidly removed, causing the smoke to
issue out of the top in puffs, as it may be seen coming out of the funnel of a railway engine. By these
means, and by varying the length of the intervals between the puffs, which had a special significance, involved messages could be sent.

\item \label{cc} Observer (Adelaide, SA)  13 Apr 1929 p. 19.  
 Codes in Smoke. 
 
The Aborigines are adept in producing ``smokes'' that convey meanings to distant observers. They have codes that are well understood over  hundreds of miles. The ``bush telegraph'' has beaten the wires sometimes with news that was important
to them. Invitations to distant tribes to corroboree or defiance to an enemy
are smoke-written in the sky. A warning
signal usually is a slender column of pale-coloured smoke. A heavy column  means ``Come at once''.

\item \label{dd} The Land (Sydney), 28 Dec 1934 p. 11. 
Smoke Telegraphy of the Aborigines. By E. S. Sorenson.

Aborigines were always on the lookout for signal smokes. 
Smoke was the Aborigines' telegraph, and they could manipulate it with wonderful skill. 

In some \textit{towris}, 
when a hunter had discovered an abundance of game he called his brothers  by sending up a  column of smoke, followed after a pause by a single cloud or puff.
A straight column without  a succeeding cloud meant that an enemy or a stranger was coming; a column
 broken at regular intervals announced the coming of friends. An enumerated series of small clouds was sent up by covering the fire with a skin or sheet of bark or a thick bush and uncovering smartly. 
Spirals were made by moving the material around the rising smoke. The corroboree signal was a row of several columns at even distances apart; and the danger signal was two columns ascending for a minute or so without a break. With accurate use of materials black puffs could be  sent aloft alongside a column of white smoke. 

My Aboriginal companion went on to a nearby hill where there was a hollow tree. 
That was his telegraph office used for communicating over long distances. He made a fire in the butt with sticks of green stuff, and a straight column of smoke went up from the top. This at intervals was varied by dabbing a sheet of bark over the fire, causing the column to break into puffs. When I remarked that the other fellow mightn't be taking notice, the operator merely grinned. His confidence was not misplaced, for the
other fellow was waiting at the crossing when I got there. 

So it was with all white travellers, smoke signals kept the Aborigines ahead, from towri to towri, apprised of their coming. 
In parts of the back country bushmen sometimes used smoke signals, but
very few could claim  knowledge of the Aborigines' code. 

\item \label{dd2} Sydney Morning Herald,  4 Apr 1942 p 7.  
Aborigines as  Japanese Foes, by Michael Terry. 

From wherever a Japanese landing may be tried a system of espionage will
radiate, by runners, smoke signals, and message stick. 
Smoke signals can get messages through the bush at extraordinary speed. Using a combination of various shapes of signals, the Aborigines ``talk'' as freely as we do by telegram.
Many times I have seen 20 to 30 smokes ascending into the sky at the
same time, and even then I have not seen an Aboriginal person or their track for days. No white man has ever made or can read the famous ``curly smoke'', which ascends spirally.  

I would suggest that the Army use trusted whites  to inform the  northern tribes of the Japanese intentions, and the consequences of defeat upon  freedom.  

\item \label{dd1} Sydney Morning Herald, 28 Sep 1946, p. 2. Mary Gilmore. Smoke Ring Code. 

Rings were a sign of urgency, puffs were the ordinary message code. Straight-up smoke meant safety or danger over. If Aborigines far away saw the signals they would send up messages so that others nearer could come at once. There was nothing mystic either in the making of the smoke or in the meaning. You learned the code as a ship's flag-officer learns the shipping code.

\item \label{ee} 
Chronicle (Adelaide)  5 Dec 1946 p. 51. 
News Carried In The Outback Smoke, by Winifred

It was very interesting to watch the Aborigines signalling. First, they build a large fire with green stuff to make smoke, then with another green bush they beat the smoke into varied shape twists and whirls which to them are definite code. As you travel through the north, you
may not see a living soul for days, yet high on the skyline spirals of smoke detail your every move.

\end{enumerate}

\end{document}